\documentclass[11pt,twoside]{article}
\usepackage{latexsym,amsmath}
\usepackage{indentfirst}
\usepackage{graphicx}

\numberwithin{equation}{section}
\newtheorem{Prop}{\bf Proposition}[section]
\newtheorem{Cor}{\bf Corollary}[section]
\newtheorem{defn}{\bf Definition}[section]
\newtheorem{Rem}{\bf Remark}[section]
\newtheorem{Ex}{\bf Example}[section]
\newtheorem{Th}{Theorem}[section]
\newtheorem{Lem}{Lemma}[section]

\begin{document}
\def \b{\Box}

\begin{center}
{\Large {\bf Study of the stability of  the fractional Stokes  system from nonlinear optics around the zero equilibrium state\\[0.1cm]
 }}

\end{center}

\begin{center}
{\bf Mihai IVAN}
\end{center}

\setcounter{page}{1}
\pagestyle{myheadings}

{\small {\bf Abstract}. The main purpose of this paper is to study the fractional-order system  with Caputo derivative associated to single Stokes pulse. The dynamic behavior for this fractional model (called the fractional Stokes system) is investigated, including: the asymptotic stability around zero equilibrium state, the stabilization problem using appropriate linear controls and the numerical integration based on fractional Euler method.
  {}

\section{Introduction}
The theory of fractional differential equations (i.e. fractional calculus) is a interesting field of mathematics based on a generalization of integer derivatives to fractional-order (\cite{nabu05, podl99}).
 The study of this field is very important and attractive because several phenomena have been proven to be better described by fractional derivatives that take into account not only the local properties,  but also global correlations of dynamical systems (\cite{gimo06, sczc19}). The fractional calculus has deep and natural connections with many fields of science and engineering (\cite{mlag, miv24, miva24}).

In the last three decades, one increasing attention has been paid to the study of the dynamic behaviors (in particular, of a series  of fractional-order differential systems associated to classical differential systems (in particular, Hamilton-Poisson systems).  These fractional models played an important role in  applied mathematics ((\cite{giva22, mivi25}), mathematical physics (\cite{giom13, migi18}), applied physics (\cite{hilf00, mivan24}),  study of biological systems (\cite{ahss07, pagci09}), chaos synchronization (\cite{igmp11, zhli05}) and so on.
   A famous Hamilton-Poisson system is the single Stokes pulse (\cite{daho90}) which comes from nonlinear optics. 

Considering the important applications of groupoids in algebra, topology, differential geometry, geometric mechanics, in applications to systems of differential equations and fractional calculus, a large number of researchers have paid great attention to studying this mathematical concept (\cite {ivan24, miva03}).

The  topological groupoids,  Lie groupoids and  Lie algebroids  have proven to be powerful tools for geometric formulation of the Hamiltonian mechanics (\cite{demi, miv03, miva05, mart}). Also, they have been used  in the investigation of many fractional dynamical systems (\cite{migo09}).

A new generalization of the notion of groupoid  is defined by Mihai Ivan in paper \cite{miva23a}, called almost groupoid. This algebraic structure can be used to define the concepts of almost topological groupoid and almost Lie groupoid.

This paper is structured as follows. A short presentation of the single Stokes pulse from nonlinear optics $ (2.9) $ is given  in Section 2. In Section 3 we define the fractional Stokes system $(3.5), $  associated to $ (2.9). $  For the system  $ (2.9) $ are proposed four types of fractional Stokes systems which are physically inequivalent. The Section 4 is dedicated to analyze of asymptotic stability of zero equilibrium state of  the fractional model $ (3.5).$  For  stabilization problem of the system $~(3.5),$ we associate the fractional Stokes  system with controls $ (4.1).$ In Theorems $ (4.1), (4.2), (4.3) $  and Corollary $ (4.1) $ are established sufficient conditions on parameters  $~k_{j},~ j=\overline{1,4}~$  to control the chaos in the fractional systems $ (4.3), (4.4), (4.6)~$ and $ (4.8).$ Using the fractional Euler's method, the numerical integration of the fractional system $(5.5) $ is presented in Section 5.\\[-0.6cm]

\section{Some theoretical aspects on the single pulse Stokes system}

We start by recalling some preliminaries about the
Hamilton-Poisson systems on ${\bf R}^{n}$. For details on
Poisson geometry and Hamiltonian dynamics see e.g. \cite{gunu93, giva13, igim11}.

 Let $C^{\infty}({\bf R}^{n}, {\bf R})$ the ring of smooth real
valued functions defined on ${\bf R}^{n}$. We shall denote a
coordinate system on ${\bf R}^{n}$ with $ x^{i} : {\bf R}^{n} \to
{\bf R},~i=\overline{1,n}.$

 A {\it Hamilton-Poisson system on ${\bf R}^{n}$} is a triple $({\bf R}^{n}, \{\cdot, \cdot \},
 H)$, where $\{\cdot, \cdot \}$ is a Poisson bracket on ${\bf
 R}^{n}$ and $H\in C^{\infty}({\bf R}^{n}, {\bf R})$ is the
 {\it Hamiltonian} ({\it energy}). Its dynamics is described by
 the following system of differential equations:
\begin{equation}
\dot{X}(t) = P(x(t))\cdot \nabla H(x(t)),\label{(1)}
\end{equation}
where $ X(t) = ( x^{1}(t), \ldots,x^{n}(t))^{T}$, $P(x(t))$ is the
structure matrix associated to Poisson bracket $\{\cdot, \cdot \}$
and $~\nabla H(x)= ( \displaystyle\frac{\partial H}{\partial
x^{1}},\ldots, \frac{\partial H}{\partial x^{n}})^{T}.$

A function  $C\in  C^{\infty}({\bf R}^{n}, {\bf R}) $ is a {\it
Casimir} of the configuration $({\bf R}^{n}, \{\cdot, \cdot
\}), $ if
\begin{equation}
P(x)\cdot \nabla C(x)=0.\label{(2)}
\end{equation}

\markboth{M. Ivan}{Study of the stability of the  fractional Stokes system from nonlinear optics ...}

We say that a dynamical system  of the form:
\begin{equation}
 \displaystyle\frac{dx^{i}(t)}{dt}= f_{i}(x^{1}(t),x^{2}(t),...,
x^{n}(t)),~~~f_{i}\in C^{\infty}({\bf R}^{n},~ {\bf
R}),~~i=\overline{1,n},  \label{(3)}
\end{equation}
 has a {\it Hamilton-Poisson
realization} (\cite{puca02}), denoted by $( {\bf R}^{n}, \Pi, H)$, if
there exist a Poisson structure $\{\cdot, \cdot\}$ on ${\bf
R}^{n}$ generated by the skew-symmetric matrix $\Pi$ and the Hamiltonian $ H $
such that the system $(2.3)$ can be written in the form $(2.1)$.

Next, we give a brief description of the single pulse Stokes system. This dynamical system is a beautiful and useful example of a Hamilton-Poisson system in nonlinear optics.

The equations of motion for the Stokes polarization parameters of a single optical beam propagating as a traveling wave in a nonlinear medium (the single pulse Stokes) from the nonlinear optics are described by using the Stokes vector $ ~u $ and a Hamiltonian function $~H.$

The Stokes vector $~u=(u^{1}, u^{2}, u^{3} )~$ is defined using the Pauli spin matrices and is called { \it the polarization parameters of the single Stokes pulse.} The vector $~u\in {\bf R}^{3} ~$ is assumed to be expressed in a linear polarization basis (\cite{daho90}). Let be $ A $ the transition matrix from the canonical basis of $ {\bf R}^{3} $  to the polarization basis. Since $ A $ is symmetric, then it one can always transform to a polarization basis in which $ A $ has the diagonal form $~ W=diag(l_{1}, l_{2}, l_{3}),$ where $~l_{i},~i=\overline{1,3}~$ are the eigenvalues of $~A.$

The Hamiltonian function H is determined by the Stokes vector $ u,$ the diagonal matrix $ W, $ and the constant vectors $ a = (a_{1}, a_{2}, a_{3}) $ and $ c =(c_{1}, c_{2}, c_{3}). $   Using the vectors $ a $ and $ c, $ we define the vector $~b= (b_{1}, b_{2}, b_{3}),$ by:\\[-0.4cm]
\begin{equation}
 b  = a + r c,~~~~~\ r= |u |.\label{(4)}
\end{equation}

The matrix $~W~$ describes the self-induced ellipse rotation. The vectors $ a $ and $ c $ describes the effects of linear and nonlinear anisotropy, respectively.

In terms of the Stokes parameters (the components of the vector $~u $), the Hamiltonian function $~H\in C^{\infty} ({\bf R}^{3}, {\bf R} ), u \rightarrow H(u), $ is defined by (\cite{daho90}) :\\[-0.2cm]
\begin{equation}
 H(u)  = u . ({\bf b} + \frac{1}{2} W . {\bf u}),~~~~ \label{(5)}
\end{equation}
where $~W = diag(l_{1}, l_{2}, l_{3}),~u = ( u^{1}~ u^{2}~u^{3 }),~{\bf b} = ( b_{1}~ b_{2}~b_{3 })^{T}~$ and
$~{\bf u} = u^{T}.$\\[-0.2cm]

 The diagonal matrix $ W $ and the choice of the vectors $ a $ and $ c, $ generates the dynamics of the Stokes vector $ u $ with the frequence $ b.$

In the coordinate system $~O u^{1} u^{2} u^{3},~$ the Hamiltonian function $ H ~ $  defined by $~(2.5),$ is written as:\\[-0.4cm]
\begin{equation}
 H(u)  = \frac{1}{2} [l_{1} (u^{1})^{2} + l_{2} (u^{2})^{2} + l_{3} (u^{3})^{2}] + b_{1} u^{1} + b_{2} u^{2}+ b_{3} u^{3}  .~~~~ \label{(6)}
\end{equation}

The dynamics of a single Stokes pulse is written as Hamilton-Poisson system. More precisely, the Stokes pulse system  is defined on the Lie-Poisson manifold $~so(3)^{*} $ (the dual of Lie algebra $ so(3) $) with  the following bracket:\\[-0.3cm]
\begin{equation}
 \{f,g\}(u) = u\cdot (\displaystyle\frac{\partial f}{\partial u}\times \displaystyle\frac{\partial g}{\partial u}), ~~~ (\forall) f, g\in C^{\infty}( {\bf R}^{3}, {\bf R}) \label{(7)}
\end{equation}
and the Hamiltonian function $~H: so(3)^{*} \cong {\bf R}^{3} \rightarrow {\bf R} ~$ given by $ (2.6).$

The dynamical system defined on  $~so(3)^{*} $ with Poisson bracket $\{., .\} $ given by $(2.7),$ enabling the equations of motion to be expressed in Hamiltonian form:\\[-0.3cm]
\begin{equation}
 \dot{u} = \{u, H \}~~ \label{(8)}
\end{equation}
where $~u\in {\bf R}^{3},~ \dot{u}= \dfrac{du}{dt}, t~$ is the time and $~H $ is the Hamiltonian function \cite{trsi87}.

Performing the calculations in the equations $~\dot{u}^{i}(t) = \{u^{i}(t), H \}~, i=\overline{1,3}~$ of the system $(2.8),$  we obtain the following differential system  on ${\bf R}^{3} $ (\cite{miva22}):\\[-0.2cm]
\begin{equation}
\left\{ \begin{array} {lcl}
 \dot{u}^{1}(t) & = & (l_{2}-l_{3})
 u^{2}(t) u^{3}(t) + b_{2}u^{3}(t)- b_{3}u^{2}(t) \\[0.1cm]
 \dot{u}^{2}(t) & = & (l_{3}-l_{1})
 u^{3}(t) u^{1}(t) + b_{3}u^{1}(t)- b_{1}u^{3}(t) \\[0.1cm]
  \dot{u}^{3}(t) & = & (l_{1}-l_{2})
 u^{1}(t) u^{2}(t) + b_{1}u^{2}(t)- b_{2}u^{1}(t), \label{(9)}
  \end{array}\right.
\end{equation}
 where the parameters  $ ~ l_{i}, b_{i}\in {\bf R} $ for $ i=\overline{1,3} $ are connected with the nature of the material and the medium (\cite{daho90, trsi87}).

The system $(2.9)$ is called the {\it single pulse Stokes  system}.
For short, the single pulse Stokes system will be called the {\it Stokes system}.

Acording to \cite{trsi87} there are six types of equations $ (2.9) $ which are physically inequivalent and which correspond to different types of optical media (\cite{puca02}).\\[-0.4cm]
\begin{Prop}
 The Stokes  system (2.9) has the  Hamilton-Poisson realization $ ({\bf R}^{3}, P, H)~$  with the Casimir  $~C \in C^{\infty}({\bf R}^{3}, {\bf R}),$  where  $ H $ is given by (2.6),\\[-0.2cm]
\begin{equation}
P(u) =  \left ( \begin{array}{ccc} 0 &
u^{3} & -u^{2}\\ \\
-u^{3} & 0 & u^{1}\\ \\
u^{2}& - u^{1} & 0
\end{array}\right )~~~\hbox{and}~~~C(u) = \displaystyle\frac{1}{2} \left (  (u^{1})^{2} + (u^{2})^{2}
+(u^{3})^{2}\right).
\label{(10)}
\end{equation}
\end{Prop}
{\it Proof.} We have $~\displaystyle\frac{\partial
H(u)}{\partial u^{1}} = l_{1} u^{1}+ b_{1},~ \displaystyle\frac{\partial
H(u)}{\partial u^{2}} = l_{2} u^{2}+ b_{2}~\displaystyle\frac{\partial
H(u)}{\partial u^{3}} = l_{3} u^{3}+ b_{3}.$ Then:\\[-0.1cm]
\[
  P(u)\cdot \nabla H(u) =
\left ( \begin{array}{ccc} 0 &
u^{3} &  - u^{2}\\[0.2cm]
-u^{3} & 0 & u^{1}\\[0.2cm]
u^{2} &  -u^{1} & 0\\[0.2cm]
\end{array}\right ) \left ( \begin{array}{c}
l_{1} u^{1} + b_{1} \\[0.2cm]
 l_{2} u^{2} + b_{2}\\[0.2cm]
  l_{3} u^{3} + b_{3}\\[0.2cm]
\end{array}\right )=\\[-0.1cm]
\]
\[
= \left (\begin{array}{c}
 (l_{2}-l_{3})
 u^{2} u^{3} + b_{2}u^{3}- b_{3}u^{2} \\[0.2cm]
(l_{3}-l_{1})
 u^{3} u^{1} + b_{3}u^{1}- b_{1}u^{3}\\[0.2cm]
(l_{1}-l_{2})
 u^{1} u^{2} + b_{1}u^{2}- b_{2}u^{1} \\[0.2cm]
\end{array}\right )=\left (\begin{array}{c}
\dot{u}^{1} \\[0.2cm]
\dot{u}^{2}\\[0.2cm]
\dot{u}^{3}\\[0.1cm]
\end{array}\right ).
\]
Hence  $ \dot{U}(t) = P(u) \cdot \nabla
H (u), $ where $\dot{u}(t) = ( \dot{u}^{1}(t),
\dot{u}^{2}(t),\dot{u}^{3}(t))^{T}$ and $(2.9)$ is a
Hamilton-Poisson system. Also, $ \nabla C(u)= ( u^{1}~~~u^{2}~~~u^{3} )^{T}.$
By a direct computation it follows $~P(u)\cdot \nabla C(u)= 0. $ Therefore,
 $ C $ is a Casimir.\hfill$\Box$\\[-0.2cm]

\section{ The fractional Stokes system}

For basic knowledge on fractional calculus, one may refer to \cite{bagr06, diet10, ivmi16, migi18, migo09,  kilb06}.

In this section we first present some theoretical aspects related to the analysis of the stability of a system of fractional differential equations around its equilibrium states.

We consider the fractional derivative operator $~ D_{t}^{q}~$ with $ q\in (0,1) $ to be {\it Caputo's
derivative}. This fractional derivative operator is often used in concrete applications.

Let $ f\in C^{\infty}(\textbf{R}) $ and $ q \in \textbf{R}, q
> 0. $ The $ q-$order Caputo differential operator
\cite{diet10}, is described by $~D_{t}^{q}f(t) = I^{m -
q}f^{(m)}(t), ~q > 0,~$ where $~f^{(m)}(t)$  represents the $
m-$order derivative of the function $ f,~m \in \textbf {N}^{\ast}$
is an integer such that $ m-1 \leq q \leq m $ and $ I^{q} $ is the
$ q-$order Riemann-Liouville integral operator \cite{kilb06}, which
is expressed by $~I^{q}f(t)
=\displaystyle\frac{1}{\Gamma(q)}\int_{0}^{t}{(t-s)^{q
-1}}f(s)ds,~q > 0,$ where $~\Gamma (.) $ is the Euler gamma  function.
  If $ q =1, $ then $ D_{t}^{q}f(t) = df/dt.$

 In the Euclidean space $ {\bf R}^{n} $ with local coordinates $~\{ x^{1}, x^{2}, \ldots, x^{n}\}, $ we consider the following system of fractional differential equations:\\[-0.3cm]
\begin{equation}
D_{t}^{q}x^{i}(t) = f_{i}(x^{1}(t), x^{2}(t), \ldots,
x^{n}(t)) ,~~ i=\overline{1,n},\label{(11)}
\end{equation}
where $q\in (0,1), f_{i}\in C^{\infty}({\bf R}^{n}, {\bf R}),
~ D_{t}^{q} x^{i}(t)$ is the Caputo fractional derivative of order $ q$ for $i=\overline{1,n}$ and $t\in [0,\tau)$ is the
time.

A point $x_{e}=(x_{e}^{1}, x_{e}^{2},\ldots, x_{e}^{n})\in {\bf
R}^{n}$ is said to be {\it equilibriun point} or {\it equilibrium state} of the system
$(3.1)$, if $~D_{t}^{q}x^{i}(t) =0$ for $i=\overline{1,n}$.

The equilibrium states of the fractional dynamical system $(3.1)$
are determined by solving the following set of equations:\\[-0.4cm]
\begin{equation}
f_{i}(x^{1}(t), x^{2}(t), \ldots, x^{n}(t)) = 0 ,~~
i=\overline{1,n}.\label{(12)}
\end{equation}
\begin{defn} {\rm (\cite{podl99}) The equilibrium point $ x_{e} $ of the system $ (3.1) $ is said to be:\\
$(i)~$ {\it (locally) stable,} if for each $\varepsilon > 0,~\exists \delta >0 $ such that\\[-0.4cm]
\begin{equation}
 \|x_{e}\| < \delta ~~ \Rightarrow ~~ \|x(t)\| < \varepsilon,~~(\forall)  t > 0.\label{(13)}
\end{equation}
$(ii)~$ {\it (locally) asymptotically stable,} if it is stable and $~\lim_{t\rightarrow \infty} \|x(t)-x_{e}\|= 0, $ where $\|\cdot\|$ is the Euclidean norm.}
\end{defn}

The Jacobian matrix associated to system $(3.1)$ is: $~~J(x)=(\displaystyle\frac{\partial f_{i}}{\partial
x^{j}}),~~~i,j=\overline{1,n}.$

\begin {Prop} {\rm (\cite{mati96})} {\it Let $ x_{e} $ be an equilibrium state of fractional differential system $(3.1)$ and
$ J(x_{e}) $ be the Jacobian matrix $J(x)$ evaluated at $ x_{e}$.

 $(i)~ x_{e}$ is locally asymptotically stable, if and only if  all eigenvalues
$ \lambda(J(x_{e})) $ of  $ J(x_{e}) $ satisfy:\\[-0.4cm]
\begin{equation}
| arg(\lambda (J(x_{e}))) | > \displaystyle\frac{q\pi}{2}.\label{(14)}
\end{equation}
 $(ii)~ x_{e} $ is locally stable, if and only if either it is asymptotically stable, or the
critical eigenvalues satisfying $~| arg(\lambda (J(x_{e}))) | = \displaystyle\frac{q \pi}{2}~$ have geometric
multiplicity one.}
\end{Prop}

Acording to Proposition 3.1, it is easy to prove the following lemma.
\begin{Lem} {\rm (\cite{giva22})} {\it Let $~x_{e}~$ be an equilibrium state of the fractional model $~(3.1)~$ and $~\lambda_{i},~i=\overline{1,n}~$ the eigenvalues of $~J(x_{e}).$

If $~\lambda_{i} < 0, $ for all $~i=\overline{1,n},~$ then  $~x_{e}~$   is asymptotically stable  $~(\forall)~q\in (0,1).$}
\end{Lem}

The Hamilton-Poisson system $(2.9)$ is modeled  by the following fractional differential equations:\\[-0.3cm]
\begin{equation}
\left\{ \begin{array} {lcl}
 D_{t}^{q}{u}^{1}(t) & = & (l_{2}-l_{3})
 u^{2}(t) u^{3}(t) + b_{2}u^{3}(t)- b_{3}u^{2}(t)  \\[0.1cm]
 D_{t}^{q}{u}^{2}(t) & = & (l_{3}-l_{1})
 u^{3}(t) u^{1}(t) + b_{3}u^{1}(t)- b_{1}u^{3}(t),~~~ q \in (0,1),\\[0.1cm]
  D_{t}^{q}{u}^{3}(t) & = & (l_{1}-l_{2})
 u^{1}(t) u^{2}(t) + b_{1}u^{2}(t)- b_{2}u^{1}(t) , \label{(15)}
  \end{array}\right.
\end{equation}
 where  $~u^{i}~$ are the Stokes polarization parameters and  $~ l_{i}, b_{i}\in {\bf R},$ for $ i=\overline{1,3}.$

The system $(3.5)$ is called the {\it fractional Stokes system} associated to $ (2.9).$
\begin{Prop}
(\cite{miva22}) The initial value problem of the fractional Stokes system $(3.5)$ has a unique solution.
\end{Prop}

As with the nonlinear dynamics generated by the Stokes pulse system $(2.9)$ there are six types of fractional equations $ (3.5) $ which are physically inequivalent and which correspond to different types of optical media.

In this section we will refer to the types of fractional Stokes systems for which the parameters $~l_{i} \in {\bf R}, i=\overline{1,3} $ meet the following condition $~~l_{1} \neq l_{2} \neq l_{3} \neq l_{1}.$

In this context there are the following four types of fractional Stokes systems:\\[0.1cm]
{\bf Type 1.} $~~b= ( 0, 0 , 0)~~~~~~~~~~~~~~~~~~~~~~~~~~~~ $ and $~~~l_{1} \neq l_{2} \neq l_{3}\neq l_{1};$\\[0.1cm]
{\bf Type 2.} $~~b= ( 0, b_{2}, 0), ~ b_{2} \neq 0~~~~~~~~~~~~~~~~~ $ and $~~~l_{1} \neq l_{2} \neq l_{3}\neq l_{1};$\\[0.1cm]
{\bf Type 3.} $~~b= ( b_{1}, 0, b_{3}), ~~ b_{1} b_{3} \neq 0~~~~~~~~~~~~$ and $~~~l_{1} \neq l_{2} \neq l_{3}\neq l_{1}; $\\[0.1cm]
{\bf Type 4.} $~~b= ( b_{1}, b_{2}, b_{3}), ~ b_{1}b_{2}b_{3} \neq 0~~~~~~~~~~ $ and $~~~l_{1} \neq l_{2} \neq l_{3} \neq l_{1}.$\\[-0.2cm]

{\bf Case 1.} The {\bf fractional Stokes  system of type 1} has the following form:\\[0.1cm]
\begin{equation}
\left\{ \begin{array} {lcl}
 D_{t}^{q}{u}^{1} & = & (l_{2}-l_{3})
 u^{2} u^{3} \\[0.1cm]
 D_{t}^{q}{u}^{2} & = & (l_{3}-l_{1})
 u^{3} u^{1},~~~~~~~~~~~~~~ q \in (0,1),\\[0.1cm]
  D_{t}^{q}{u}^{3} & = & (l_{1}-l_{2})
 u^{1} u^{2}, \label{(16)}
  \end{array}\right.
\end{equation}
 where  $~ l_{i}\in {\bf R},~i=\overline{1,3} $ such that $~~~l_{1} \neq l_{2} \neq l_{3}\neq l_{1}.$

An easy computation shows us that the equilibrium states of system (3.6) are:\\
$~e_{0}=(0, 0, 0),~~e_{11}^{m}= (m, 0, 0),~~e_{12}^{m}= (0, m, 0),~~e_{13}^{m}= (0, 0, m),~~~m\in {\bf R}^{\ast}.$\\[-0.2cm]

{\bf Case 2.} The {\bf fractional Stokes  system of type 2} has the following form:\\[0.1cm]
\begin{equation}
\left\{ \begin{array} {lcl}
 D_{t}^{q}{u}^{1} & = & (l_{2}-l_{3})
 u^{2} u^{3} + b_{2}u^{3}  \\[0.1cm]
 D_{t}^{q}{u}^{2} & = & (l_{3}-l_{1})
 u^{3} u^{1} ,~~~~~~~~~~~~~~~~~~~~~~~~~~~ q\in (0,1)\\[0.1cm]
  D_{t}^{q}{u}^{3} & = & (l_{1}-l_{2})
 u^{1} u^{2} - b_{2}u^{1} , \label{(17)}
  \end{array}\right.
\end{equation}
 where  $~ b_{2}\in {\bf R}^{\ast} $ and  $ l_{i}\in {\bf R}, i=\overline{1,3} $  such that $~ l_{1} \neq l_{2} \neq l_{3} \neq l_{1}.$

An easy computation shows us that the equilibrium states of system (3.7) are:\\
$~e_{0}=(0, 0, 0),~~e_{21}^{m}= (m, \frac{b_{2}}{l_{1} - l_{2}}, 0),~~e_{3}^{m}= (0, m, 0),~~e_{23}^{m}= (0, \frac{b_{2}}{l_{3} - l_{2}}, m),~~~m\in {\bf R}^{\ast}.$\\[-0.3cm]

{\bf Case 3.} The {\bf fractional Stokes  system of type 3} has the following form:\\[0.1cm]
\begin{equation}
\left\{ \begin{array} {lcl}
 D_{t}^{q}{u}^{1} & = & (l_{2}-l_{3})
 u^{2} u^{3} - b_{3}u^{2}  \\[0.1cm]
 D_{t}^{q}{u}^{2} & = & (l_{3}-l_{1})
 u^{3} u^{1} + b_{3}u^{1}- b_{1}u^{3},~~~~~~~~~~~~~ q \in (0,1),\\[0.1cm]
  D_{t}^{q}{u}^{3} & = & (l_{1}-l_{2})
 u^{1} u^{2} + b_{1}u^{2}, \label{(18)}
  \end{array}\right.
\end{equation}
where  $~ b_{1}, b_{3}\in {\bf R}^{\ast} $ and  $ l_{i}\in {\bf R}, i=\overline{1,3} $  such that $~ l_{1} \neq l_{2} \neq l_{3} \neq l_{1}.$

An easy computation shows us that the equilibrium states of system (3.8) are:\\
$~e_{0}=(0, 0, 0),~~e_{32}^{m}= (\frac{b_{1}}{l_{2} - l_{1}}, m, \frac{b_{3}}{l_{2} - l_{3}} )~$ for $~m\in {\bf R}~ $ and\\
$e_{33}^{m}= (\frac{b_{1}m}{b_{3}+ (l_{3} - l_{1})m}, 0, m),~$ for $~m\in {\bf R}\setminus \{\frac{b_{3}}{l_{1} - l_{3}}\}.$\\[-0.3cm]

{\bf Case 4.} The {\bf fractional Stokes  system of type 4} has the following form:\\[0.1cm]
\begin{equation}
\left\{ \begin{array} {lcl}
 D_{t}^{q}{u}^{1} & = & (l_{2}-l_{3})
 u^{2} u^{3} + b_{2}u^{3}- b_{3}u^{2}  \\[0.1cm]
 D_{t}^{q}{u}^{2} & = & (l_{3}-l_{1})
 u^{3} u^{1} + b_{3}u^{1}- b_{1}u^{3},~~~~~~~~~~~~~~~~ q \in (0,1),\\[0.1cm]
  D_{t}^{q}{u}^{3} & = & (l_{1}-l_{2})
 u^{1} u^{2} + b_{1}u^{2}- b_{2}u^{1}, \label{(19)}
  \end{array}\right.
\end{equation}
where  $~ b_{1}, b_{2},  b_{3}\in {\bf R}^{\ast} $ and  $ l_{i}\in {\bf R}, i=\overline{1,3} $  such that $~ l_{1} \neq l_{2} \neq l_{3} \neq l_{1}.$

An easy computation shows us that the equilibrium states of system (3.9) are:\\
$~e_{0}=(0, 0, 0)~~$ and\\
$e_{41}^{m}= (\frac{b_{1}m}{b_{3}+ (l_{3} - l_{1})m}, \frac {b_{2}m}{b_{3}+ (l_{3} - l_{2})m}, m),~$ for $~m\in {\bf R}\setminus \{\frac{b_{3}}{l_{1} - l_{3}}, \frac{b_{3}}{l_{2} - l_{3}} \}.$

It is easy to see that  $~e_{0}=(0, 0, 0)~$ is an equilibrium state for all four types of fractional Stokes system $ (3.5).$\\[-0.3cm]
\begin{Rem}
{\rm  The fractional Stokes system of {\bf type $ 2 $} was studied in \cite{miva22}, including: the asymptotic stability of its equilibrium states and the stabilization problem using suitable linear controls.}
\end{Rem}
\begin{Rem}
{\rm  If in the fractional model $ (3.5)$ we take $ q=1, $ then one obtains the system for integer-order derivative which  corresponds to case $ 4 $ of the dynamics $ (2.9).$ For  this
 dynamical system, the nonlinear stability and the problem of existence of periodic solutions are studied, see (\cite{puca02}).} \\[-0.6cm]
\end{Rem}

\section{ Asymptotic stability of the fractional Stokes system $ (3.5) $ }

Let us we present the study of asymptotic stability of the fractional system $(3.5)$ around the zero equilibrium state. Finally, we will discuss
how to stabilize the unstable zero equilibrium states of the system $(3.5)$  (in all the four cases) via Caputo fractional  derivative. For this study we apply the Matignon's test (Proposition 3.1).\\[0.2cm]
{\bf 4.1. Stability analysis of the fractional Stokes system  $(3.5)$ around zero equilibrium state}

In this subsection we present the study  of asymptotic stability for the equilibrium state $~e_{0}~$ of the fractional system $ (3.5). $

The Jacobian matrix associated to system $(3.5)$ is:
\[
J(u, b) = \left ( \begin{array}{ccc}
0 & (l_{2}-l_{3}) u^{3} -b_{3}    &  (l_{2}-l_{3}) u^{2} + b_{2} \\
  (l_{3}-l_{1}) u^{3} +b_{3}  & 0 & (l_{3}-l_{1}) u^{1} -b_{1}\\
  (l_{1}-l_{2}) u^{2}- b_{2}   &  (l_{1}-l_{2}) u^{1} + b_{1}  & 0 \\
\end{array}\right ).\\[-0.1cm]
\]
\begin{Prop}
The equilibrium state $ e_{0}= (0, 0, 0 )~$ of the fractional system (3.5)
is unstable $ (\forall) q \in (0,1).$
\end{Prop}
{\it Proof.} The characteristic polynomial of the matrix $~
J(e_{0}, b) =\left (\begin{array}{ccc}
  0 & -b_{3} &  b_{2} \\
  b_{3}  & 0 & -b_{1}\\
   - b_{2}   &  b_{1}  & 0 \\
\end{array}\right ) $
is $~p_{J(e_{0}, b)}(\lambda) = \det ( J(e_{0}, b) -
\lambda I) = - \lambda (\lambda ^{2} +  b_{1}^{2}+ b_{2}^{2} + b_{3}^{2}).~$

 The equation $~p_{J(e_{0}, b)}(\lambda) = 0 $ has the root $ \lambda_{1} = 0 .$  Since $~arg(\lambda_{1}) =0 < \frac{q \pi}{2}$
for all $ q\in (0,1),$  by Proposition 3.1 follows that $~e_{0}~$ is unstable for all $ q\in (0,1).$ \hfill$\Box$\\[0.2cm]
{\bf 4.2.  Controllability of chaotic behaviors of the fractional Stokes system  $(3.5)$ }

In this subsection we will apply the general  method for to control the stability of the fractional system $ (3.5),~$  see \cite{giva22}.

Let  $ x_{e}$ be an unstable equilibrium point of the fractional Stokes system $~(3.5).~$ We associate to $(3.5)$ a new fractional-order system with controls and given by:\\[-0.2cm]
\begin{equation}
\left\{ \begin{array} {lcl}
 D_{t}^{q}{u}^{1}(t) & = & (l_{2}-l_{3})
 u^{2}(t) u^{3}(t) + b_{2}u^{3}(t)- b_{3}u^{2}(t) + \varphi_{1}(t) \\[0.1cm]
 D_{t}^{q}{u}^{2}(t) & = & (l_{3}-l_{1})
 u^{3}(t) u^{1}(t) + b_{3}u^{1}(t)- b_{1}u^{3}(t) + \varphi_{2}(t),~~~~~~~ q \in (0,1),\\[0.1cm]
  D_{t}^{q}{u}^{3}(t) & = & (l_{1}-l_{2})
 u^{1}(t) u^{2}(t) + b_{1}u^{2}(t)- b_{2}u^{1}(t) + \varphi_{3}(t) , \label{(20)}
  \end{array}\right.
\end{equation}
where $~\varphi_{i}(t), i=\overline{1,3}$ are control functions.

If one selects the control functions $~\varphi_{i}, ~i=\overline{1,3} $  which then make the eigenvalues of the linearized equation of $(4.1)$ satisfy one of the
conditions from Proposition 3.1, then the trajectories of $ (4.1) $ asymptotically approaches the unstable equilibrium state $x_{e}$ in the sense that
$\lim_{t\rightarrow \infty} \|x(t)-x_{e}\|= 0.$ In the case when $~x_{e}~$ is unstable, then fractional model $ (4.1) $  may exhibit chaotic behavior.

In the following we will choose the control functions $~\varphi_{i}, ~i=\overline{1,3} $ in each of the four types of fractional system (4.1) to control the chaotic behavior of
the respective system around the zero equilibrium state.\\[0.1cm]
{\bf $\bullet$}~{\it Controllability of chaotic behaviors of the fractional model $(3.5)$ of types {\bf 1} and {\bf 2}.}

If in the fractional model $(4.1) $ we take $~\varphi_{i}(t), i=\overline{1,3}, $  given by:\\[-0.4cm]
 \begin{equation}
\varphi_{1}(t) = k_{1} u^{1}(t),~~~ \varphi_{2}(t) = k_{2} u^{2}(t),~~~ \varphi_{3}(t) = k_{1} u^{3}(t), ~~~ k_{1}, k_{2} \in {\bf R}^{*}.\label{(21)}
\end{equation}
we obtain the controlled fractional Stokes systems associated to fractional systems of  $(3.6) $ and $ (3.7), $ respectively.

With the control functions $(4.2) $ and $ b_{1} = b_{2}= b_{3} = 0,~$ the system $(4.1) $ becomes:\\[-0.2cm]
\begin{equation}
\left\{ \begin{array} {lcl}
 D_{t}^{q}{u}^{1} & = & (l_{2}-l_{3})u^{2} u^{3} + k_{1} u^{1} \\[0.1cm]
 D_{t}^{q}{u}^{2} & = & (l_{3}-l_{1})
 u^{3} u^{1} + k_{2} u^{2},~~~~~~~~~~~~~ q \in (0,1),\\[0.1cm]
  D_{t}^{q}{u}^{3} & = & (l_{1}-l_{2})
 u^{1} u^{2} + k_{1} u^{3}, \label{(22)}
  \end{array}\right.
\end{equation}
 where   $~ l_{i}\in {\bf R},~i=\overline{1,3} $ such that $~~~l_{1} \neq l_{2} \neq l_{3}\neq l_{1}~$ and $~ k_{1}, k_{2} \in {\bf R}^{\ast}~$ are control parameters.

The fractional system $ (4.3) $ is called the {\it controlled fractional Stokes system associated  to
$~(3.6) $ at  $~e_{0}.$}

With the control functions $(4.2) $ and $ b_{1} = b_{3} = 0,~$ the system $(4.1) $ becomes:\\[-0.2cm]
\begin{equation}
\left\{ \begin{array} {lcl}
 D_{t}^{q}{u}^{1} & = & (l_{2}-l_{3})u^{2} u^{3} + b_{2}u^{3} + k_{1} u^{1} \\[0.1cm]
 D_{t}^{q}{u}^{2} & = & (l_{3}-l_{1})
 u^{3} u^{1} + k_{2} u^{2},~~~~~~~~~~~~~ q \in (0,1),\\[0.1cm]
  D_{t}^{q}{u}^{3} & = & (l_{1}-l_{2})
 u^{1} u^{2}- b_{2}u^{1} + k_{1} u^{3}, \label{(23)}
  \end{array}\right.
\end{equation}
 where $~b_{2}\in {\bf R}^{\ast},~ l_{i}\in {\bf R},~i=\overline{1,3} $ such that $~~~l_{1} \neq l_{2} \neq l_{3}\neq l_{1}~$ and $~ k_{1}, k_{2} \in {\bf R}^{*}~$ are control parameters.

The fractional system $ (4.4) $ is called the {\it controlled fractional Stokes system associated  to
$~(3.7) $ at  $~e_{0}.$}

The Jacobian matrix of the fractional model $(4.4)$ with the controls $ k_{1}, k_{2} $ is\\[-0.2cm]
\[
J(u, b_{2}, k_{1}, k_{2}) = \left ( \begin{array}{ccc}
k_{1} & (l_{2}-l_{3}) u^{3}    &  (l_{2}-l_{3}) u^{2} + b_{2} \\
  (l_{3}-l_{1}) u^{3}  & k_{2} & (l_{3}-l_{1}) u^{1}\\
  (l_{1}-l_{2}) u^{2}- b_{2}   &  (l_{1}-l_{2}) u^{1}  & k_{1} \\
\end{array}\right ).\\[-0.1cm]
\]
\begin{Th}
Let  be the controlled fractional Stokes  system $(4.4) $ with the controls $ k_{1}, k_{2}\in {\bf R^{\ast}}. $\\
$(i)~$ If  $~k_{1} < 0, k_{2} < 0, $ then $ e_{0}= (0, 0, 0)$ is asymptotically stable $~(\forall)~q\in (0,1).$\\
$(ii)~$ If  $~k_{1} > 0, k_{2} < 0 $  and $~q_{0}=\displaystyle\frac{2}{\pi}
\arctan\displaystyle\frac{|b_{2}|}{ k_{1}}, ~$ then:\\
$(1)~ e_{0} $ is asymptotically stable $~(\forall)~q \in ( 0, q_{0}),$
and it is stable for $~q = q_{0};$\\
$(2)~ e_{0} $ is unstable $~(\forall)~q \in ( q_{0}, 1).$\\
$(iii)~$ If  $~k_{2} > 0 $ and $~k_{1}\in {\bf R}^{\ast}, $ then $ e_{0}$ is unstable $~(\forall)~q\in (0,1).$
\end{Th}
{\it Proof.} We have $~J_{02}:= J(e_{0}, b_{2} , k_{1}, k_{2}) =\left (\begin{array}{ccc}
  k_{1} & 0 &  b_{2} \\
  0  & k_{2} & 0\\
   - b_{2}   &  0  & k_{1} \\
\end{array}\right ).~$
Its characteristic polynomial is $~p_{02}(\lambda):= p_{J_{02}}(\lambda) = \det ( J_{02} -
\lambda I) = - (\lambda - k_{2}) [(\lambda - k_{1}) ^{2} + b_{2}^{2}].$
The roots of the equation $ p_{02}(\lambda)= 0~$ are $\lambda_{1}=
k_{2},~\lambda_{2,3}= k_{1} \pm i b_{2}.$\\
$(i)~$ We suppose that $~k_{1} < 0~$ and $~k_{2} < 0.$ In
this case  we have $ Re (\lambda_{i})<
0 $ for $ i=\overline{1,3}. $ Since $|arg(\lambda_{i})| =\pi > \displaystyle\frac{q \pi}{2},
i=\overline{1,3} $ for all $ q\in (0, 1)$, by Proposition 3.1(i),
it implies that $e_{0} $ is asymptotically stable for all $ q\in (0,1).$\\
$(ii)~$ We suppose that $~k_{1} > 0, k_{2} < 0. $  In this case  we have $ \lambda_{1} < 0 $ and $ Re (\lambda_{2,3})>
0. $ Applying Proposition 3.1(i), $ e_{0} $ is asymptotically stable, for $~0 < q < q_{0},$ where $~ q_{0}=\displaystyle\frac{2}{\pi}
\arctan\displaystyle\frac{|b_{2}|}{k_{1}}.$ If $ q = q_{0},$ then $ e_{0} $ is stable. For $~ q_{0} < q < 1,
~e_{0} $ is unstable. Hence, the assertion $(ii)$ holds.\\
$(iii)~$ We suppose that $~k_{2} > 0 ~$ and $~k_{1}\in {\bf R}^{\ast}. $  In this case,
$~ J_{02} $ has at least a positive eigenvalue. For instance, since $ \lambda_{1} > 0 $ we have $~|arg(\lambda_{1}| = 0< \displaystyle\frac{q \pi}{2}~$  for all $~q\in (0,1), $ it implies that  $ e_{0} $ is unstable $ (\forall) k_{1} \in {\bf R}^{\ast}.$  Hence, $(iii) $ holds. \hfill$\Box$\\[-0.4cm]
\begin{Cor}
Let  be the controlled fractional Stokes  system $(4.3) $ with the controls $ k_{1}, k_{2}\in {\bf R^{\ast}}~ $ and $~q\in (0,1).$\\
$(i)~$ If  $~k_{1} < 0, k_{2} < 0, $ then $ e_{0}= (0, 0, 0)$ is asymptotically stable.\\
$(ii)~$ If  $~k_{1} > 0, k_{2} < 0 ,$  then $~e_{0} $ is unstable\\
$(iii)~$ If  $~k_{2} > 0 $ and $~k_{1}\in {\bf R}^{\ast}, $ then $ e_{0} $ is unstable.
\end{Cor}
{\it Proof.} The characteristic polynomial of Jacobian matrix $~J_{01}: = J(e_{0}, k_{1}, k_{2})~$ is\\
$~p_{01}(\lambda):= p_{J_{01}}(\lambda) = - (\lambda - k_{2}) (\lambda - k_{1}) ^{2}.$
The roots of the equation $ p_{01}(\lambda)= 0~$ are $\lambda_{1}=
k_{2},~\lambda_{2,3}= k_{1}.$ Since $~b_{2}=0, $ it implies  $~q_{0}=0. $  Now, the statements $~(i), (ii)~$ and $(iii)~$ of this corollary follow immediately from Theorem 4.1. \hfill$\Box$\\[-0.4cm]

The eigenvalues of the fractional model $(4.3), $ {\it before control} are $~\lambda_{1,2,3} =0 ~$ and {\it after control} are $~ \lambda_{1}= k_{2},~ \lambda_{2,3}= k_{1}. $
\begin{Ex}
{\rm  Let be the  controlled fractional Stokes system $(4.3). $
We select $~ l_{1} =-1, l_{2} = 0.58, l_{3} = -0.08. $ \\
$~(i) $ Chosing $ k_{1} =-0.64, k_{2} = -0.43 ~$ and  according to Corollary $4.1(i) $ it follows that  $~ e_{0} = (0, 0,0)~$ is asymptotically stable for $~ q=0.68.$\\
$~(ii) $  Choose $~k_{1}= 0. 32 $ and  $ k_{2} = -0.06. $  By Corollary $ 4.1(ii), e_{0}= ( 0, 0, 0) $ is unstable for $~ q=0.98.$}\\[-0.2cm]
\end{Ex}
In Table 1 we give a set of values for  $~l_{i}, k_{1}, k_{2}, i=\overline{1,3}, ~$  and the corresponding eigenvalues of controlled fractional Stokes  system $(4.3), $ before control and after control.
{\small
$$\begin{array}{|l|c|c||c|c|c|} \hline
   l_{1}, l_{2}, l_{3} & \lambda_{i}  & Stability & controls  & \lambda_{i}  &  Stability \cr
& before~ control & q\in (0,1) &  k_{1}, k_{2} & after~ control & q\in (0,1)\cr \hline

l_{1} =1, l_{2} = 0.55, &  &  &  k_{1}=-1.35  & \lambda_{1}=-0.65 &  asym.  \cr
l_{3} =- 0.45    & \lambda_{1,2,3}= 0 & unstable & k_{2}=-0.65 & \lambda_{2,3}= -1.35 &  stable  \cr \hline

l_{1} =1, l_{2} = 0.55, &  &  & k_{1}=-2  & \lambda_{1}=1.2  &    \cr
l_{3} =- 0.45  & \lambda_{1,2,3}= 0 & unstable & k_{2}=1.2 & \lambda_{2,3}=-2 &  unstable  \cr \hline
\end{array} $$\\[-0.2 cm]

{\bf Table 1.} {\it The controls $~k_{1}, k_{2},~$ equilibrium state $~e_{0}~$ and corresponding eigenvalues}.}\\[-0.2cm]

The eigenvalues of the fractional model $(4.4), $  {\it before control} are $~\lambda_{1} =0 ,~ \lambda_{2,3}=\pm i b_{2}~$ and {\it after control} are $~ \lambda_{1}= k_{2},~ \lambda_{2,3}= k_{1}\pm i b_{2}. $
\begin{Ex}
{\rm  Let be the  controlled fractional Stokes  system $(4.4). $
We select $~ l_{1} =-0.9, l_{2} = 0.15, l_{3} = 1.05 ~$ and  $~b_{2}=-0.81. $ \\
$~(i) $ Chosing $ k_{1} =-0.13 ~$ and $~k_{2} = -0.24, ~$ by Theorem $ 4.1 (i) $  it follows that $~ e_{1} = (0, 0,0)~$ is asymptotically stable for $~ q=0.64.$\\
$~(ii) $ For $~k_{1}= 0.81, k_{2} = -0.02,$  it follows that $~ q_{0} = 0,5. $ Then, according to Theorem 4.1(ii),~$  e_{0}= ( 0, 0, 0) $ is  asymptotically stable for $~ q=0.48 $ and unstable for $~ q=0.81.$}\\[-0.4cm]
\end{Ex}
In Table 2 we give a set of values for  $~l_{i}, b_{2}, k_{1}, k_{2}, i=\overline{1,3}, ~$  and the corresponding eigenvalues of controlled fractional Stokes  system $(4.4), $ before control and after control.
{\small
$$\begin{array}{|l|c|c||c|c|c|} \hline
   l_{1}, l_{2}, l_{3}, b_{2} & \lambda_{i}  & Stability & controls  & \lambda_{i}  &  Stability \cr
& before~ control & q\in (0,1) &  k_{1}, k_{2} & after~ control & q\in (0,1)\cr \hline

l_{1} =1, l_{2} =- 0.55, & \lambda_{1} = 0 &  &  k_{1}=-1.8  & \lambda_{1}=-0.35 &  asym.  \cr
l_{3} =- 0.45, b_{2}=1    & \lambda_{2,3}=\pm i  & unstable & k_{2}=-0.35 & \lambda_{2,3}= -1.8\pm i &  stable  \cr \hline

l_{1} =1, l_{2} = 0.55, & \lambda_{1} = 0 &  & k_{1}=0.5773  & \lambda_{1}=-2  & asym.   \cr
l_{3} =- 0.45,  & \lambda_{2,3}=\pm i  & unstable & k_{2}=-2 & \lambda_{2,3}=0.5773\pm i &  stable  \cr
b_{2}=-1 &  & & & q_{0}=0.66 &  q\in (0, 0.66)  \cr \hline

l_{1} =1, l_{2} = 0.55, & \lambda_{1} = 0 &  & k_{1}=0.5773  & \lambda_{1}=-2  &    \cr
l_{3} =- 0.45,  & \lambda_{2,3}=\pm i  & unstable & k_{2}=-2 & \lambda_{2,3}=0.5773\pm i &  unstable  \cr
b_{2}=-1 &  & & & q_{0}=0.66 &  q\in (0.66,1)  \cr \hline

l_{1} =1, l_{2} =- 0.55, & \lambda_{1} = 0 &  &  k_{1}=-0.12  & \lambda_{1}=0.72 &   \cr
l_{3} =- 0.45, b_{2}=-0.24   & \lambda_{2,3}=\pm i  & unstable & k_{2}= 0.72 & \lambda_{2,3}= -0.12 \pm i &  unstable  \cr \hline
\end{array} $$\\[-0.3 cm]

{\bf Table 2.} {\it The controls $~k_{1}, k_{2},~$ equilibrium state $~e_{0}~$ and corresponding eigenvalues}.}\\[0.1cm]

{\bf $\bullet$} {\it Controllability of chaotic behaviors of the fractional model  $(3.5)$ of type {\bf 3}}

If in the fractional model $(3.5) $ we take $~ \varphi_{i}(t), i=\overline{1,3}, $  given by:\\[-0.2cm]
 \begin{equation}
\varphi_{1}(t) = k_{3} u^{1}(t),~~~ \varphi_{2}(t) = 0,~~~ \varphi_{3}(t) = k_{3} u^{3}(t), ~~~ k_{3}\in {\bf R}^{\ast}.\label{(24)}
\end{equation}
we obtain the controlled fractional Stokes systems associated to  $(3.8).$

With the control functions $(4.6)~$ and $~b_{2}= 0, $ the system $(4.1) $ becomes:\\[-0.2cm]
\begin{equation}
\left\{ \begin{array} {lcl}
 D_{t}^{q}{u}^{1} & = & (l_{2}-l_{3})u^{2} u^{3} - b_{3} u^{2} + k_{3} u^{1} \\[0.1cm]
 D_{t}^{q}{u}^{2} & = & (l_{3}-l_{1})
 u^{3} u^{1} + b_{3} u^{1} - b_{1} u^{3},~~~~~~~~~~~~~ q \in (0,1),\\[0.1cm]
  D_{t}^{q}{u}^{3} & = & (l_{1}-l_{2})
 u^{1} u^{2} +  b_{1} u^{2} + k_{3} u^{3}, \label{(25)}
  \end{array}\right.
\end{equation}
 where $~ b_{1}, b_{3}\in {\bf R}^{\ast},~ l_{i}\in {\bf R},~i=\overline{1,3} $ such that $~~~l_{1} \neq l_{2} \neq l_{3}\neq l_{1}~$ and $~ k_{3}\in {\bf R}^{\ast}~$ is a control parameter.

The fractional system $ (4.6) $ is called the {\it controlled fractional Stokes system associated  to
$~(3.8) $ at  $~e_{0}.$}

The Jacobian matrix of the fractional model $(4.6)$ with the control $ k_{3} $ is\\[-0.2cm]
\[
J(u, b_{1}, b_{3}, k_{3}) = \left ( \begin{array}{ccc}
k_{3} & (l_{2}-l_{3}) u^{3} - b_{3}    &  (l_{2}-l_{3}) u^{2} \\
  (l_{3}-l_{1}) u^{3} +b_{3} & 0 & (l_{3}-l_{1}) u^{1} - b_{1}\\
  (l_{1}-l_{2}) u^{2}  &  (l_{1}-l_{2}) u^{1} + b_{1}  & k_{3} \\
\end{array}\right ).\\[-0.1cm]
\]
\begin{Th}
Let  be the controled fractional Stokes  system $(4.6) $ with the control $ k_{3}\in {\bf R^{\ast}},~ \beta_{03} =\sqrt{b_{1}^{2}+ b_{3}^{2}}~$
 and  $ \Delta_{03} = k_{3}^{2} -4 \beta_{03}^{2}.~$\\
 ${\bf 1}.~$ Let $~\Delta_{03} < 0 ~$ and $~q\in (0, 1).$\\
$(i)~$ If $~k_{3}<0, $ then $ e_{0}$ is asymptotically stable $~(\forall)~k_{3}\in (-2 \beta_{03}, 0).$\\
$(ii)~$ If  $~k_{3} > 0, $ then $ e_{0}~ $ is unstable $~(\forall)~k_{3} \in (0, 2 \beta_{03}).$\\
${\bf 2}.~$ Let $~\Delta_{03} > 0 ~$ and $~q\in (0, 1).$\\
$(i)~$ If $~k_{3}<0, $ then $ e_{0}$ is asymptotically stable $~(\forall)~k_{3} \in (-\infty, -2 \beta_{03}].$\\
$(ii)~$ If  $~k_{3} > 0, $ then $ e_{0}~ $ is unstable $~(\forall)~k_{3} \in [ 2\beta_{03}, \infty).$
\end{Th}
{\it Proof.} We have $~J_{03}:= J(e_{0}, b_{1}, b_{3}, k_{3}) =\left (\begin{array}{ccc}
  k_{3} & -b_{3} &  0\\
  b_{3}  & 0 & -b_{1}\\
  0 & b_{1}   & k_{3} \\
\end{array}\right ).~$ Its characteristic polynomial is $~p_{03}(\lambda):= p_{J_{03}}(\lambda) = \det ( J_{03} -
\lambda I) = - (\lambda - k_{3}) (\lambda^{2} - k_{3}\lambda  + b_{1}^{2} + b_{3}^{2} ).$
The roots of the equation $ p_{03}(\lambda)= 0~$ are $\lambda_{1}= k_{3},~~\lambda_{2,3}= \frac{k_{3}\pm \sqrt{\Delta_{03}}}{2},~$ where $ \Delta_{03} = k_{3}^{2} -4 (b_{1}^{2} + b_{3}^{2}).~$ If we denote
$~\beta_{03} = \sqrt{b_{1}^{2} + b_{3}^{2}},~$ then $ \Delta_{03} = k_{3}^{2} -4 \beta_{03}^{2}.~$\\
{\bf Case 1.} $~\Delta_{03} < 0~$  and $~q\in (0,1).~$  We have $~\Delta_{03} < 0~$  if and only if $~k_{3}\in (-2 \beta_{03}, 2 \beta_{03}).~$  Then $~ \lambda_{1}= k_{3}, \lambda_{2,3}= \frac{k_{3}\pm i \sqrt{-\Delta_{03}}}{2}.~$

$(i)~$ We suppose  $~ k_{3} < 0.~$ In this case  we have $ \lambda_{1} < 0 $ and $ Re (\lambda_{2,3})<
0.$ Since $|arg(\lambda_{i})| =\pi > \displaystyle\frac{q \pi}{2},~i=\overline{1,3}, $ by Proposition 3.1(i), it implies that $~e_{0} $ is  asymptotically stable for all $~k_{3}\in (-2 \beta_{03}, 0).~$\\
$(ii)~$ We suppose  $~ k_{3} > 0.~$  Since $|arg(\lambda_{1})| =0 < \displaystyle\frac{q \pi}{2}, $ by Proposition 3.1(i), it implies that $~e_{0} $ is  unstable  for all $~k_{3}\in (0, 2 \beta_{03}).~$
Therefore, the assertion $(ii) $ holds.

{\bf Case 2.} $~\Delta_{03} \geq 0~$  and $~q\in (0,1).~$  We have $~\Delta_{03} \geq 0~$  if and only if $~k_{3}\in (-\infty, -2 \beta_{03}]\cup [ 2\beta_{03}, \infty).~$  Then
 $~ \lambda_{1}= k_{3},~\lambda_{2,3}=\frac{k_{3}\pm \sqrt{ \Delta_{03}}}{2},~$

$(i)~$ Let  $~ k_{3} < 0.~$ Then  $~\lambda_{1} < 0. $  We have $~ \lambda_{2} + \lambda_{3}= k_{3}<0 ~$ and $~ \lambda_{2}\lambda_{3}= \beta_{03}^{2}>0.~$ It follows that  $ ~\lambda_{2}<0~$
and $~\lambda_{3}<0. $  Applying Corollary 4.1, it implies that $~e_{0} $ is  asymptotically stable for all $~k_{3}\in (-\infty, -2 \beta_{03}].$

$(ii)~$ Let $~k_{3} > 0.~$  Then  $~\lambda_{1} >0.~$ Then, $~J(e_{0}, k_{3})~$ has at least a positive eigenvalue and so $ e_{0}~$ is unstable $~(\forall)~k_{3}\in [ 2\beta_{03}, \infty).~$  Hence, the assertion $(ii) $ holds.\hfill$\Box$\\[-0.2cm]

The eigenvalues of the fractional model $(4.6), $  with $~\Delta_{03} = k_{3}^{2}- 4 \beta_{03}^{2},~ $  are:\\
{\it before control}: $~\lambda_{1} =0 ,~ \lambda_{2,3}=\pm i \beta_{03}~$  and {\it after control}:\\
{\bf Case 1.~} for $~\Delta_{03}<0 ~~\Rightarrow ~~ \lambda_{1}= k_{3}, ~\lambda_{2,3}=\frac {k_{3}\pm i \sqrt{-\Delta_{03}}}{2};$ \\
{\bf Case 2.~} for $~\Delta_{03}>0 ~~\Rightarrow ~~ \lambda_{1}= k_{3},~\lambda_{2,3}=\frac {k_{3}\pm \sqrt{\Delta_{03}}}{2}. $
\begin{Ex}
{\rm  Let be the  controlled fractional Stokes  system $(4.6). $
We select $~ l_{1} =-0.6, l_{2} = 0.35, l_{3} = 1.01 ~$ and  $~b_{1}=0.3 , b_{3} =0.4. $ We have $~\beta_{03}= 0.5 $ \\
$~(i) $ Choose $ k_{3} =-0.6 .~$  Then $~\Delta_{03}= -0.64. ~$  By Theorem $ 4. 1(i) $  it follows that $~ e_{0} = (0, 0,0)~$ is asymptotically stable for $~ q=0.85.$\\
$~(ii) $ If $~k_{3}= -2,~$ then $~\Delta_{03}= 3. $  by Theorem $ 4. 1(ii) $  it follows that $~ e_{0} = (0, 0,0)~$ is asymptotically stable for $~ q=0.73.$\\
$~(iii) $ If $~k_{3}= 3,~$ then $~\Delta_{03}= 2. $  By Theorem $ 4. 2(ii) $  it follows that $~ e_{0} = (0, 0,0)~$ is unstable for $~ q=0.43.$}\\[-0.4cm]
\end{Ex}

In Table 3 we give a set of values for  $~l_{i}, b_{1}, b_{3}, k_{3}, i=\overline{1,3}, ~$  and the corresponding eigenvalues of controlled fractional model $(4.6), $ before control and after control.

{\small
$$\begin{array}{|l|c|c||c|c|c|} \hline
   l_{1}, l_{2}, l_{3}, b_{1}, b_{3} & \lambda_{i}  & Stability & control  & \lambda_{i}  &  Stability \cr
& before~ control & q\in (0,1) &  k_{3} & after~ control & q\in (0,1)\cr \hline

l_{1} =1, l_{2} =- 0.2, & \lambda_{1} = 0 &  &  & \lambda_{1}=-0.8 &  asym.  \cr
l_{3} =- 0.42, b_{1}=1.2, & \lambda_{2,3}=\pm 1.3416 i  & unstable & k_{3}=-0.8 & \lambda_{2,3}= &  stable  \cr
b_{3}=0.6 &  & &  &  -0.4\pm 1.7581 i &   \cr \hline

l_{1} =0.6, l_{2} =- 0.12, & \lambda_{1} = 0 &  &  & \lambda_{1}=1.2 &   \cr
l_{3} =1.1, b_{1}=1.2, & \lambda_{2,3}=\pm 1.3416 i  & unstable & k_{3}=1.2 & \lambda_{2,3}= &  unstable  \cr
b_{3}=0.6 &  & &  &  0.6\pm 1.7581 i &   \cr \hline

l_{1} =0.8, l_{2} =- 0.35, & \lambda_{1} = 0 &  &  & \lambda_{1}=-1.2 &  asym.  \cr
l_{3} =0.28, b_{1}=0.3, & \lambda_{2,3}=\pm 0.5 i  & unstable & k_{3}=-1.2 & \lambda_{2}=-0.264 &  stable  \cr
b_{3}=0.4 &  & &  & \lambda_{3}= -0.9316 &   \cr \hline

l_{1} =0.8, l_{2} =- 0.35, & \lambda_{1} = 0 &  &  & \lambda_{1}=2 &   \cr
l_{3}= 0.28, b_{1}=0.3, & \lambda_{2,3}=\pm 0.5 i  & unstable & k_{3}= 2 & \lambda_{2}= 1,866&  unstable  \cr
b_{3}=0.4 &  & &  &  \lambda_{3}= 0.134 &   \cr \hline
\end{array} $$\\[-0.2 cm]

{\bf Table 3.} {\it The control $~k_{3},~$ equilibrium state $~e_{0}~$ and corresponding eigenvalues}.}\\[0.1cm]

{\bf $\bullet$} {\it Controllability of chaotic behaviors of the fractional model $(3.5)$ of type {\bf 4}}

If in the fractional model $(4.1) $ we take  $~ \varphi_{i}(t), i=\overline{1,3}, $  given by:\\[-0.2cm]
 \begin{equation}
\varphi_{1}(t) = k_{4} u^{1}(t),~~~ \varphi_{2}(t) = k_{4} u^{2}(t),~~~ \varphi_{3}(t) = k_{4} u^{3}(t), ~~~ k_{3}\in {\bf R}^{\ast}.\label{(26)}
\end{equation}
we obtain the controlled fractional Stokes systems associated to $(3.9) .$

With the control functions $(4.7),~$  the system $(4.1) $ becomes:\\[-0.2cm]
\begin{equation}
\left\{ \begin{array} {lcl}
 D_{t}^{q}{u}^{1} & = & (l_{2}-l_{3})u^{2} u^{3} + b_{2} u^{3} - b_{3} u^{2} + k_{4} u^{1} \\[0.1cm]
 D_{t}^{q}{u}^{2} & = & (l_{3}-l_{1})
 u^{3} u^{1} + b_{3} u^{1} - b_{1} u^{3} + k_{4} u^{2},~~~~~~~~~~~~~ q \in (0,1),\\[0.1cm]
  D_{t}^{q}{u}^{3} & = & (l_{1}-l_{2})
 u^{1} u^{2} +  b_{1} u^{2} - b_{2} u^{1} + k_{4} u^{3}, \label{(27)}
  \end{array}\right.
\end{equation}
 where $~ b_{1}, b_{2}, b_{3}\in {\bf R}^{\ast},~ l_{i}\in {\bf R},~i=\overline{1,3} $ such that $~~~l_{1} \neq l_{2} \neq l_{3}\neq l_{1}~$ and $~ k_{4}\in {\bf R}^{\ast}~$ is a control parameter.

The fractional system $ (4.8) $ is called the {\it controlled fractional Stokes system associated  to
$~(3.9) $ at  $~e_{0}.$}

The Jacobian matrix of the fractional model $(4.8)$ with the control $ k_{4} $ is\\[-0.2cm]
\[
J(u, b, k_{4}) = \left ( \begin{array}{ccc}
k_{4} & (l_{2}-l_{3}) u^{3} - b_{3}    &  (l_{2}-l_{3}) u^{2} + b_{2} \\
  (l_{3}-l_{1}) u^{3} + b_{3} & k_{4} & (l_{3}-l_{1}) u^{1} - b_{1}\\
  (l_{1}-l_{2}) u^{2} - b_{2}  &  (l_{1}-l_{2}) u^{1} + b_{1}  & k_{4} \\
\end{array}\right ).\\[-0.1cm]
\]
\begin{Th}
Let  be the controlled fractional Stokes  system $(4.8) $ with the control $ k_{4}\in {\bf R^{\ast}}~$ and $~\beta_{04} =\sqrt{b_{1}^{2}+ b_{2}^{2} + b_{3}^{2}}.$\\
$(i)~$ If $~k_{4}<0, $ then $~e_{0}~$ is asymptotically stable $~(\forall)~q\in (0,1).$\\
$(ii)~$ If  $~k_{4} > 0, $ then $ e_{0}~ $ is unstable $~(\forall)~q\in (0,1).$
\end{Th}
{\it Proof.} We have $~J_{04}:= J(e_{0}, b,  k_{4}) =\left (\begin{array}{ccc}
  k_{4} & -b_{3} & b_{2}\\
  b_{3}  & k_{4} & -b_{1}\\
  -b_{2} & b_{1}   & k_{4} \\
\end{array}\right ).~$ Its characteristic polynomial is $~p_{04}(\lambda):= p_{J_{04}}(\lambda) = \det ( J(e_{0}, k_{4}) -
\lambda I) = - (\lambda - k_{4}) (\lambda^{2} + b_{1}^{2} + b_{2}^{2} + b_{3}^{2} ).$
The roots of the equation $ p_{04}(\lambda)= 0~$ are $\lambda_{1}= k_{4},~\lambda_{2,3}= k_{4}\pm i\beta_{04},~$   where
$~\beta_{04} = \sqrt{b_{1}^{2} + b_{2}^{2} + b_{3}^{2}}.$\\
$(i)~$ Let  $~ k_{4} < 0.~$ In this case  we have $ \lambda_{1} < 0 $ and $ Re (\lambda_{2,3})<
0. $ By Proposition 3.1(i), it implies that $~e_{0} $ is  asymptotically stable. Therefore, the assertion $(i) $ holds.\\
$(ii)~$ Let $~k_{4} > 0. ~$  In this case, $~ J(e_{0}, k_{4}) $ has  a positive eigenvalue. It implies that  $ e_{0} $ is unstable. Hence, assertion $(ii) $ holds.
\hfill$\Box$\\[-0.2cm]

The eigenvalues of the fractional model $(4.8),~$  {\it before control} are $~\lambda_{1} =0 ,~ \lambda_{2,3}=\pm i \beta_{04}~$  and {\it after control} are
$\lambda_{1}= k_{4}, ~\lambda_{2,3}= k_{4}\pm i \beta_{04}.$
\begin{Ex}
{\rm  Let be the  controlled fractional Stokes  system $(4.8). $
We select $~ l_{1} =-0.16, l_{2} = 0.56, l_{3} = 1.04 ~$ and  $~b_{1}=0.3 , b_{2}= 0.6, b_{3} =0.6. $ We have $~\beta_{04}= 0.9 $ \\
$~(i) $ Choose $ k_{4} =-1.4 .~$  By Theorem $4. 3.(i) $  it follows that $~ e_{0} = (0, 0,0)~$ is asymptotically stable for $~ q=0.48.$\\
$~(ii) $ If $~k_{4}= 2.1~$ then, by Theorem $ 4.3.(ii), $  it follows that $~ e_{0} = (0, 0,0)~$ is unstable for $~ q=0.37.$}\\[-0.4cm]
\end{Ex}

In Table 4 we give a set of values for  $~l_{i}, b_{i}, k_{4}, i=\overline{1,3}, ~$  and the corresponding eigenvalues of controlled fractional model $(4.8), $ before control and after control.
{\small
$$\begin{array}{|l|c|c||c|c|c|} \hline
   l_{1}, l_{2}, l_{3},  & \lambda_{i}  & Stability & control  & \lambda_{i}  &  Stability \cr
b_{1}, b_{2}, b_{3} & before~ control & q\in (0,1) &  k_{4} & after~ control & q\in (0,1)\cr \hline

l_{1} =1, l_{2} =- 0.5, & \lambda_{1} = 0 &  &  & \lambda_{1}=-0.83 &  asym.  \cr
l_{3} =- 0.64, b_{1}=0.1, & \lambda_{2,3}=\pm 0.3 i  & unstable & k_{4}= -0.83 & \lambda_{2,3}= &  stable  \cr
b_{2}=0.2 , b_{3}= 0.2 &  & &  & -0.83 \pm 0.3 i  &   \cr \hline

l_{1} =0.6, l_{2} =- 0.12, & \lambda_{1} = 0 &  &  & \lambda_{1}=1.32 &   \cr
l_{3} =1.1, b_{1}=1.2, & \lambda_{2,3}=\pm 0.3 i  & unstable & k_{4}= 1.32 & \lambda_{2,3}= &  unstable  \cr
b_{2}=0.2 , b_{3}=0.2 &  & &  &  1.32\pm 0.3 i &   \cr \hline

\end{array} $$\\[-0.2 cm]

{\bf Table 4.} {\it The control $~k_{4},~$ equilibrium states $~e_{0}~$ and corresponding eigenvalues}.}\\[-0.5cm]

\section{Numerical integration of the fractional Stokes system with controls $(4.1)$}

In this section we start with some mathematical preliminaries of the fractional Euler's method for solving initial value problem for fractional differential equations.

Consider the following general form of the initial value problem (IVP) with Caputo derivative \cite{odmo08}:\\[-0.4cm]
\begin{equation}
D_{t}^{q} y(t) = f(t,y(t)),~~~ y(0)=y_{0},~~~t\in I=[0,T],~T>0  \label{(28)}
\end{equation}
where $~y: I \rightarrow {\bf R}^{n},~f: {\bf R}^{n} \rightarrow {\bf R}^{n}~$ is a continuous nonlinear function and
$ q\in (0,1), $ represents the order of the derivative.

The right-hand side of the IVP $ (5.1) $ in considered examples are Lipschitz functions and the numerical method used in this works to integrate system $ (5.1) $ is the {\it Fractional Euler's method}.

Since $ f $ is assumed to be continuous function, every solution of the initial value problem given by $ (5.1) $ is also a solution of the following {\it Volterra fractional integral equation}:\\[-0.4cm]
\begin{equation}
y(t) = y(0) +~ I_{t}^{q}f(t,y(t)), \label{(29)}
\end{equation}
where $ I_{t}^{q} $ is the $ q-$order Riemann-Liouville integral operator, which is expressed by:\\[-0.4cm]
\begin{equation}
I_{t}^{q}f(t) =\displaystyle\frac{1}{\Gamma(q)}\int_{0}^{t}{(t-s)^{q
-1}}f(s,y(s))ds,~q > 0,~~~s\in [0, T].  \label{(30)}
\end{equation}
Moreover, every solution of $ (5.2) $ is a solution of the (IVP) $ (5.1).$

To integrate the fractional equation $( 5.1),$ means to find the solution of $ (5.2) $ over the interval $~[0,T]. $ In this context, a set of points $~(t_{j}, y(t_{j})) $ are produced which are used as approximated values. In order to achieve this approximation, the interval $ [0,T] $  is partitioned into $ n $ subintervals $~[t_{j}, t_{j+1} ] $  each equal width $~ h =\frac{T}{n}, ~ t_{j} = j h $  for $ j = 0,1,..., n. $

For the fractional-order $ q $ and $~ j = 0, 1, 2, ... ,$ it computes an approximation denoted as $~y_{j+1} ~$ for $~ y ( t_{j+1}),~ j= 0, 1, ... .$

The general formula of the fractional Euler's method for to compute the elements $~y_{j},$ is\\[-0.4cm]
\begin{equation}
y_{j+1} = y_{j} + \displaystyle\frac{h^{q}}{\Gamma(q+1)} f( t_{j}, y(t_{j})),~~~~~ t_{j+1} = t_{j} + h,~~~ j=0, 1, ..., n. \label{(31)}
\end{equation}

For the numerical integration of the fractional Stokes  system $ (4.1),$ we apply the fractional Euler method (FEM). For this, consider the following fractional differential equations\\[-0.2cm]
\begin{equation}
\left\{\begin{array}{lcl}
 D_{t}^{q} u^{i}(t) & = &
F_{i}(u^{1}(t), u^{2}(t), u^{3}(t)),~~~ i=\overline{1,3}, ~~t\in
(t_{0}, \tau),~q \in (0,1)\\
u(t_{0}) &=& (u^{1}(t_{0}), u^{2}(t_{0}), u^{3}(t_{0}))
\end{array}\right.\label{(32)}\\[-0.2cm]
\end{equation}
where\\[-0.4cm]
\begin{equation}
\left\{\begin{array}{lcl} F_{1}(u(t)) & = & (l_{2}-l_{3})
 u^{2}(t) u^{3}(t) + b_{2}u^{3}(t) - b_{3}u^{2}(t) + d_{1} u^{1}(t),\\
F_{2}(u(t)) & = &  (l_{3}-l_{1})
 u^{3}(t) u^{1}(t) + b_{3}u^{1}(t)- b_{1}u^{3}(t) + d_{2} u^{2}(t),\\
F_{3}(u(t)) & = & (l_{1}-l_{2})
 u^{1}(t) u^{2}(t) + b_{1}u^{2}(t)- b_{2}u^{1}(t) + d_{3} u^{3}(t),\\
\end{array}\right.\label{(33)}
\end{equation}
 where  $~b_{i}, d_{i}\in {\bf R}~$ and  $ l_{i}\in {\bf R}, i=\overline{1,3} $  such that $~ l_{1} \neq l_{2} \neq l_{3} \neq l_{1}.$\\

 Since the functions $ F_{i}(u(t)), i=\overline{1,3} $ are continuous, the initial value problem $(5.5)$ is equivalent
 to system of Volterra integral equations, which is given as follows:\\[-0.2cm]
\begin{equation}
u^{i}(t)~=~ u^{i}(0)  + ~I_{t}^{q} F_{i}(u^{1}(t), u^{2}(t), u^{3}(t)),
~~~~~i=\overline{1,3}.\label{(34)}
\end{equation}

The system $(5.7)$ is called the {\it Volterra integral equations
associated to fractional Stokes system} $(4.1)$.

The problem for solving the system $(5.5)$ is reduced to one of solving a sequence of systems of fractional equations in
increasing dimension on successive intervals  $[j, (j+1)]$.

For the numerical integration of the system $(5.6)$ one can use the fractional Euler method (the formula $(5.4)$\emph{}), which is expressed as follows:\\[-0.4cm]
\begin{equation}
u^{i}(j+1)=u^{i}(j)+ \frac{h^{q}}{\Gamma (q+1)} F_{i}(u^{1}(j), u^{2}(j),
u^{3}(j)),~~~ i=\overline{1,3},\label{(35)}
\end{equation}
where $ j=0,1,2,...,N,  h=\displaystyle\frac{T}{N}, T>0, N>0.~$

More precisely, the numerical integration of the fractional system $(5.8)$ is given by:\\[-0.4cm]
\begin{equation}
\left \{ \begin{array}{ll} u^{1}(j+1) &= u^{1}(j)+
h^{q}~\displaystyle\frac{1}{\Gamma(q+1)}(\alpha
 u^{2}(j) u^{3}(j) + b_{2} u^{3}(j) - b_{3} u^{2}(j) + d_{1} u^{1}(j))\\[0.3cm]
u^{2}(j+1) &= u^{2}(j)+
h^{q}~\displaystyle\frac{1}{\Gamma(q+1)}( \beta u^{1}(j) u^{3}(j) + b_{3} u^{1}(j) - b_{1} u^{3}(j)+ d_{2} u^{2}(j))\\[0.3cm]
 u^{3}(j+1) &= u^{3}(j) +
h^{q}~\displaystyle\frac{1}{\Gamma(q+1)}(\gamma
 u^{1}(j) u^{2}(j) + b_{1} u^{2}(j) - b_{2} u^{1}(j) +  d_{3} u^{3}(j))\\[0.3cm]
u^{i}(0)&= u_{e}^{i}+\varepsilon,~~~i=\overline{1,3},\\
\end{array}\right. \label{(36)}
\end{equation}
where $~\alpha := l_{2} -l_{3},~ ~\beta := l_{3} -l_{1},~ ~\gamma := l_{1} -l_{2}.$

Using \cite{diet10, odmo08}, we have that the numerical algorithm given by
$(5.9)$ is convergent.\\[-0.2cm]

If in the relations $ (5.9) $ we replace the parameters $~b_{i}, d_{i}, i=\overline{1,3}~$ with the values corresponding to the four types of fractional systems studied,
we wiil obtain the numerical integration of the fractional systems with controls $~(4.3), (4.4), (4.6)~$ and $~(4.8).~$ More precisely, for:\\
$-~~~ b_{1}=b_{2}=b_{3}=0~$ and $~d_{1}=k_{1}, d_{2}=k_{2}, d_{3}= k_{1}~$ we obtain the numerical integration  $~(5.9).1 $ of the fractional model with controls $~(4.3);$\\
$-~~~ b_{1}=b_{3}=0~$ and $~d_{1}=k_{1}, d_{2}=k_{2}, d_{3}= k_{1}~$ we obtain the numerical integration $~ (5.9).2 $ of the fractional model with controls $~(4.4);$\\
$-~~~ b_{2}=0~$ and $~d_{1}=k_{3}, d_{2}=0, d_{3}= k_{3}~$ we obtain the numerical integration  $~ (5.9).3 $ of the fractional model with controls $~(4.6);$\\
$-~~~ d_{1}=k_{4}, d_{2}=k_{4}, d_{3}= k_{4}~$ we obtain the numerical integration $~ (5.9).4 $ of the fractional model with controls $~(4.8).$\\[-0.5cm]
\begin{Ex}
{\rm  Let us we present the numerical integration of the fractional Stokes  system with controls $ (4.8) $  which has considered in Example $~4.4(i).$ For this we apply the algorithm $ (5.9).4 $ and software Maple. Then, in $(5.9).4~$ we take: $~ l_{1} =-0.16, l_{2} = 0.56, l_{3} = 1.04, b_{1}= 0.3, b_{2}= 0.6,  b_{3}= 0.6~$  and  $~k_{4} = -1.4. $ It is known that the equilibrium state $~ e_{0} = (0, 0,0)~$ is asymptotically stable for $~ q=0.48.$

For the numerical simulation of solutions of the above fractional model we use the rutine {\it Maple. fract-Stokes-system-with-controls, type 4}, denoted by [sp-fr.Stokes-syst 4]. Applying this program for
$~h = 0.01, \varepsilon= 0.01, u^{1}(0) = \varepsilon , u^{2}(0)= \varepsilon, u^{3}(0)=\varepsilon, N = 100, t = 102 $ one obtain the orbits $~(n, u^{1}(n)), (n, u^{2}(n)), (n, u^{3}(n))~$ and  $( (u^{1}(n), u^{2}(n), u^{3}(n) ) $ , for $ q = 0.48.$}
\end{Ex}

Finally, we present the rutine  [sp-fr.Stokes-syst 4]:\\[0.2cm]
$ \# $ Fractional equations associated to Stokes system for q=0.48\\[0.1cm]
Du1/dt=(l2-l3)*u2*u3 + b2*u3 -b3*u2 + k4* u1; Du2/dt=(l3-l1)*u1*u3 + b3*u1-b1*u3+ k4* u2; Du3/dt=(l1-l2)*u1*u2 + b1*u2 - b2*u1 + k4* u3;\\[0.2cm]
{\small
$>$ with (plots):

$>$ l1:=-0.16; l2:=0.56; l3:=1.04; alpha:=l2-l3; beta:=l3-l1; gamma:=l1-l2;
 b1:=0.3; b2:=0.6; b3:=0.6; k4:= -1.4; q:=0.48; u1e:=0.; u2e:=0.; u3e:=0.;\\[0.1cm]
$>$ with (stats):

$>$ h:=0.01; epsilon:=0.01;
 n:=100:t:=n+2; u1:= array (0 .. n): u2:= array (0 .. n): u3:= array (0 .. n): u1[0]:=epsilon + u1e; u2[0]:=epsilon + u2e; u3[0]:=epsilon + u3e;\\[0.1cm]
$>$ for ~j~ from ~1~  by ~1~ to ~n~  do

$>$ u1[j]:= u1[j-1] + h $\wedge$ q *(alpha* u2[j-1]*u3[j-1] + b2*u3[j-1] - b3*u2[j-1] + k4* u1[j-1])/GAMMA(q+1);

 u2[j]:= u2[j-1] + h $\wedge$ q *(beta* u1[j-1]*u3[j-1] + + b3*u1[j-1] - b1*u3[j-1] +  k4* u2[j-1])/GAMMA(q+1);

 u3[j]:= u3[j-1] + h $\wedge$ q *(gamma* u1[j-1]*u2[j-1] + b1*u2[j-1] - b2*u1[j-1] + k4* u3[j-1])/GAMMA(q+1);\\[0.1cm]
od:

$>$ plot (seq([j,u1[j]], j = 0 .. n), style = point, symbol = point, scaling = UNCONSTRAINED);

plot (seq([j,u2[j]], j = 0 .. n), style = point, symbol = point, scaling = UNCONSTRAINED);

plot (seq([j,u3[j]], j = 0 .. n), style = point, symbol = point, scaling = UNCONSTRAINED);

pointplot 3d ( {seq([u1[j], [u2[j], [u3[j]],  j = 0 .. n)}, style = point, symbol = point, scaling = UNCONSTRAINED, color = red);}\\[-0.5cm]
\begin{Rem}
{\rm  Appyling  $~(5.9)~$ and the package Maple for the numerical simulation of solutions of fractional models $~(4.3), (4.4), (4.6) $ and $ (4.8)~$ for each set of values for parameters $~ l_{i}, b_{i}~ $ for $~i=\overline{1,3}~$ and $~k_{j},~j=\overline{1,4}~$ given in the Tables 1-4, it will be found that the results obtained are valid.}\\[-0.4cm]
\end{Rem}

{\bf Conclusions.} This paper presents the fractional Stokes  system  $(3.5) $ associated to  system $(2.9).$  The fractional Stokes system $ (3.5) $ was studied from fractional differential equations theory point of view: asymptotic stability, determining of sufficient conditions on parameters $ k_{j},~j=\overline{1,4} $ to control the chaos in the controlled fractional system (in the four cases) associated to $ (3.5) $ and  numerical integration of the controlled fractional model $ (4.1).$ The obtained results in the paper was illustrated with examples for their feasibility. The study of chaotic fractional systems has applications in theory of chaos synchronization and secure communications. In this context,  by
 choosing the right control functions $ \varphi_{i}, i=\overline{1,3}~ $ in the controlled fractional Stokes system  $~(4.1),~$ this work offers a series of chaotic and  non-chaotic  fractional differential systems.\\

Author's adress\\[-0.2cm]

West University of Timi\c soara. Geometry and Topology Seminar.\\
Teacher Training Department. Timi\c soara. Romania.\\
E-mail: mihai.ivan@e-uvt.ro\\
\end{document}